\documentclass[11pt]{article}
\usepackage{amscd, amsmath, amssymb, epic}

\title{Hilbert's fourteenth problem over finite fields,
and a conjecture on the cone of curves}
\author{Burt Totaro}
\date{  }

\def\Z{\text{\bf Z}}
\def\Q{\text{\bf Q}}
\def\R{\text{\bf R}}

\def\P{\text{\bf P}}
\def\F{\text{\bf F}}

\def\Pic{\text{Pic}}
\def\arrow{\rightarrow}

\def\qed{\ QED }

\def\Pic{\text{Pic}}
\def\Hom{\text{Hom}}

\def\half{\frac{1}{2}}
\def\diag{\text{diag}}

\def\Gal{\text{Gal}}
\def\Aut{\text{Aut}}
\def\Bir{\text{Bir}}
\def\PsAut{\text{PsAut}}
\def\Proj{\text{Proj\:}}
\def\Int{\text{Int\,}}
\def\NE{\overline{NE}}
\def\A{\overline{A}}

\begin{document}
\maketitle

\newtheorem{theorem}{Theorem}[section]
\newtheorem{corollary}[theorem]{Corollary}
\newtheorem{lemma}[theorem]{Lemma}
\newtheorem{conjecture}[theorem]{Conjecture}

Hilbert's fourteenth problem asks whether the ring of invariants
of any representation of a linear algebraic group
is finitely generated over the base field.
Nagata gave the first counterexample,
using a representation of $(G_a)^{13}$, where
$G_a$ denotes the additive group \cite{NagataICM}.
In his example, the representation is defined over a field
of large transcendence degree over the prime field (of any
characteristic). Mukai simplified Nagata's construction,
showing that there are representations of $(G_a)^3$ over the
complex numbers whose ring of invariants is not finitely
generated \cite{Mukai}. 

Mukai relates the problem
of finite generation of rings of invariants to a natural question
in algebraic geometry: when is the total coordinate ring 
of a projective variety finitely
generated? The deepest known result is the
Birkar-Cascini-Hacon-M\textsuperscript{c}Kernan theorem that the total
coordinate ring of a Fano variety in characteristic zero
is finitely generated \cite[Cor.\ 1.3.1]{BCHM}.
Nagata and Mukai give counterexamples
to finite generation for certain varieties (blow-ups of projective
space) just outside the realm of Fano varieties. The precise border
between finite and infinite generation remains to be understood.

Mukai's construction yields representations whose coefficients
are ``general'' complex numbers. In this paper, we show
that Mukai's three best examples
can all be realized over finite
fields and over the rational numbers,
in fact with simple explicit coefficients.
We give examples over all fields, including the field of order 2.
These are the first published counterexamples to Hilbert's
fourteenth problem over finite fields. (In retrospect,
some examples over finite fields
can be constructed using the work of Manin \cite{Manin},
Tate \cite{Tate}, or Shioda \cite{Shioda}
on elliptic surfaces.) For some classes
of examples, we characterize exactly when finite generation holds
and when it does not.

\begin{theorem}
\label{main}
Let $k$ be any field.
Then there are
linear representations over $k$ of $(G_a)^3$ on $A^{18}$,
of $(G_a)^4$ on $A^{16}$, and of $(G_a)^6$ on $A^{18}$ whose rings
of invariants are not finitely generated. The representations are defined
explicitly.
\end{theorem}

Here is the geometric idea. The basic ingredient
of Nagata's examples (especially as simplified by Steinberg \cite{Steinberg})
is the existence of non-torsion line bundles of degree zero on an elliptic
curve. There are no such line bundles on curves over
finite fields. But the generic fiber of an elliptic fibration
over a finite field can have infinite Mordell-Weil group,
and that turns out to be enough. Precisely, we encounter known elliptic
fibrations of blow-ups of $\P^2$ and $\P^3$, and what seems to be
a new fibration of a blow-up of $\P^5$ by abelian surfaces
(section \ref{abelian}). It would be interesting to describe the geometry
of these rational abelian fibrations of $\P^5$
in more detail. For example, are these fibrations completely integrable
systems with respect to some Poisson structure?

Our main results (Theorems \ref{sharp2} and \ref{sharp3})
relate finite generation
of the total coordinate ring, in some situations,
to finiteness of a certain Mordell-Weil
group. We conclude with a more general conjecture,
saying that the cone of curves of any variety with semi-ample
anticanonical bundle
is controlled by a group,
which may be infinite.
This would follow from a generalization of the
conjectures of Kawamata and Morrison
on Calabi-Yau fiber spaces to allow klt pairs (Conjecture \ref{domain}).
We prove the conjecture
for smooth projective surfaces with semi-ample anticanonical bundle,
the new case being that of rational elliptic surfaces (Theorem
\ref{rational}).

Thanks to Igor Dolgachev, Shigeru Mukai,
and Greg Sankaran for their comments.

\section{Another example}

The smallest known representation of an algebraic group for which
finite generation fails is Freudenburg's
11-dimensional representation of a unipotent group $(G_a)^4\rtimes G_a$
over the rational numbers \cite{Freudenburg},
based on an example by Kuroda \cite{Kuroda}. It would be interesting
to know whether there are such low-dimensional examples over finite fields.

There are broader forms of Hilbert's fourteenth problem,
for example about actions of algebraic groups on arbitrary affine varieties.
Since even the most specific form of the problem, about
linear representations, has a negative answer,
we focus on that case.

\section{Mukai's method}

In this section, we summarize Mukai's geometric approach to
producing counterexamples to Hilbert's fourteenth problem.

We use the following result by Mukai \cite{Mukai}. (He works
over the complex numbers, but his proof is elementary 
and works over any field.) 

\begin{theorem}
\label{mukai}
For any $n\geq r\geq 3$,
let $X$ be the blow-up of projective space $\P^{r-1}$ at $n$
distinct rational points $p_1,\ldots,p_n$, not contained
in a hyperplane, over a field $k$.
Let $G\cong (G_a)^{n-r}$ be the subgroup of $(G_a)^n$
which is the kernel of a linear map $A^n\arrow A^r$ corresponding
to the points $p_i$. Let $(G_a)^n$ act on $V=A^{2n}$ by
$$(t_1,\ldots,t_n)(x_1,\ldots,x_n,y_1,\ldots,y_n)=(x_1,\ldots,x_n,
y_1+t_1x_1,\ldots,y_n+t_nx_n).$$
Then the ring of invariants $O(V)^G$ is isomorphic to the total
coordinate ring of $X$,
$$TC(X):=\oplus_{a,b_1,\ldots,b_n\in\Z} H^0(X,aH-b_1E_1-\cdots
-b_nE_n)\cong \oplus_{L\in \Pic X}H^0(X,L),$$
where $H$ is the pullback of the hyperplane line bundle on $\P^{r-1}$
and $E_1,\ldots,E_n$ are the exceptional divisors in $X$.
\end{theorem}

Therefore, to give examples of rings of invariants which are not
finitely generated, it suffices to exhibit arrangements of points
on projective space such that the total coordinate ring of
the blow-up $X$ is not finitely generated. 

Define a {\it pseudo-isomorphism }between smooth
projective varieties to be a birational map which is an isomorphism
outside subsets of codimension at least 2. This notion can be useful
for singular varieties,
but our applications only involve smooth varieties.
(Example: a pseudo-isomorphism
between smooth projective surfaces
is an isomorphism.) A pseudo-isomorphism $X\dashrightarrow
Y$ induces an isomorphism $\Pic(X)\arrow \Pic(Y)$ by taking the
proper transform of divisors.
Define a {\it $(-1)$-divisor }$D$ on a projective variety $X$
to be the proper transform under some pseudo-isomorphism $X\dashrightarrow
X'$ of the exceptional divisor for a morphism $X'\arrow Y$
which is the blow-up
of a smooth point. (For example, a $(-1)$-curve on a surface
is simply a curve isomorphic to $\P^1$ with self-intersection $-1$,
by Castelnuovo.)
Mukai observes that every $(-1)$-divisor is indecomposable
in the monoid of effective line bundles on $X$ (line bundles $L$
with $H^0(X,L)\neq 0$). Also, two different
$(-1)$-divisors represent different elements of $\Pic (X)$. Therefore,
if a variety $X$ contains infinitely
many $(-1)$-divisors, then the monoid of effective line bundles on
$X$ is not finitely generated, and so the total coordinate ring of $X$
is not finitely generated \cite[Lemma 3]{Mukai}. 

Lemma \ref{cremona} will give a way to ensure that a blow-up $X$ of projective
space contains infinitely many $(-1)$-divisors. To state it,
define a birational
map $\Psi: \P^{r-1}\dashrightarrow \P^{r-1}$ by
$$[x_1,\ldots,x_r]\mapsto \bigg[\frac{1}{x_1},\ldots,\frac{1}{x_r}\bigg].$$
It contracts the $r$ coordinate hyperplanes to the $r$ coordinate
points. A birational map between two projective spaces
which is projectively equivalent to $\Psi$
is called a {\it standard Cremona transformation}. 
We note that $\Psi$ lifts to a pseudo-isomorphism
from the blow-up of $\P^{r-1}$ at the $r$ coordinate points to
itself. In modern terms, this pseudo-isomorphism can be described
as a composite of several flops.

We say that an arrangement of $n$ points in $\P^{r-1}$, $n\geq r\geq 3$,
is in {\it linear
general position }if no $r$ of the points are contained in a hyperplane.
(In particular, the points are all distinct.)
Given $n$ points in linear general position, we can perform
the standard Cremona transformation on any $r$ of the $n$ points;
this gives a different arrangement of $n$ points in 
a projective space.
They need not be in linear general position. We say that an arrangement
of $n$ points in $\P^{r-1}$ is in {\it Cremona general position }if
they are in linear general position and this remains true
after any finite sequence of standard Cremona transformations
on $r$-tuples of the points.

The key point of Mukai's method is:

\begin{lemma}
\label{cremona}
 Let $p_1,\ldots, p_n$, $n\geq r\geq 3$, be points in projective
space $\P^{r-1}$ over a field $k$ which are in Cremona general position.
If
$$\frac{1}{2}+\frac{1}{r}+\frac{1}{n-r}\leq 1,$$
then the blow-up $X$ of $\P^{r-1}$ at $p_1,\ldots,p_n$ contains
infinitely many $(-1)$-divisors. Therefore the total coordinate
ring of $X$ is infinitely generated, and the corresponding representation
of $(G_a)^{n-r}$ of dimension $2n$ over $k$ has infinitely generated
ring of invariants.
\end{lemma}

For clarity, we recall the proof of Lemma \ref{cremona}.
Clearly $X$ contains
$n$ $(-1)$-divisors over the $n$ points $p_1,\ldots,p_n$. But since
$p_1,\ldots,p_n$ are in Cremona general position, we can perform
Cremona transformations on $r$-tuples of points in $p_1,\ldots,p_n$
any number of times. This shows that $X$ is pseudo-isomorphic
to a blow-up of $\P^{r-1}$ in many other ways, and so we find other
$(-1)$-divisors on $X$. The inequality on $n$
ensures, by a purely combinatorial argument, that the resulting
$(-1)$-divisors on $X$ have arbitrarily large degrees (when projected
down to divisors on $\P^{r-1}$), and so there are infinitely many of them.
To describe this argument in more detail, define
a symmetric bilinear form on $\Pic(X)$ by $H^2=r-2$, $E_i^2=-1$,
$H\cdot E_j=0$,
and $E_i\cdot E_j=0$ for $i\neq j$. We have
a canonical identification of $\Pic(X)$ with $\Z^{n+1}=\Z H\oplus
\Z E_1\oplus \cdots \oplus \Z E_n$ for every blow-up $X$ of
an $(r-1)$-dimensional projective space
at an ordered set of $n$ points. So the standard Cremona transformation
at the points $p_1,\ldots,p_r$ 
determines a well-defined automorphism of $\Z^{n+1}$, which we compute
to be 
the reflection $s_n$ orthogonal to $H-E_1-\cdots-E_r$.
Switching points $p_i$ and $p_{i+1}$
acts on $\Z^{n+1}$ by the reflection $s_i$ orthogonal to $E_i-E_{i+1}$
for $1\leq i\leq {n-1}$. We compute (it is also clear geometrically)
that the Cremona action on $\Pic(X)=\Z^{n+1}$ fixes the anticanonical class
of $X$, $-K_X=rH-(r-2)E_1-\cdots-(r-2)E_n$.
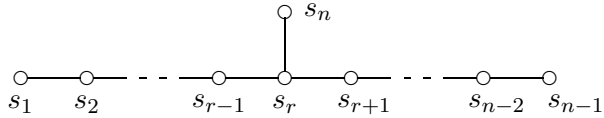
\begin{figure}[ht] \label{dynkin1}
\hspace{2cm}
\begin{picture}(300,60)
\put(50,20){\line(1,0){20}}
\put(75,20){\line(1,0){10}}
\dashline{3}(85,20)(110,20)
\put(110,20){\line(1,0){10}}
\put(125,20){\line(1,0){20}}
\put(150,20){\line(1,0){20}}
\put(175,20){\line(1,0){10}}
\dashline{3}(185,20)(210,20)
\put(210,20){\line(1,0){10}}
\put(225,20){\line(1,0){20}}
\put(147.5,22.5){\line(0,1){20}}
\put(47.5,20){\circle{5}}
\put(72.5,20){\circle{5}}
\put(122.5,20){\circle{5}}
\put(147.5,20){\circle{5}}
\put(172.5,20){\circle{5}}
\put(222.5,20){\circle{5}}
\put(247.5,20){\circle{5}}
\put(147.5,45){\circle{5}}
\put(47.5,30){\makebox(0,0){}}
\put(197.5,30){\makebox(0,0){}}
\put(47.5,10){\makebox(0,0){$s_{1}$}}
\put(72.5,10){\makebox(0,0){$s_{2}$}}
\put(122.5,10){\makebox(0,0){$s_{r-1}$}}
\put(147.5,10){\makebox(0,0){$s_{r}$}}
\put(177.5,10){\makebox(0,0){$s_{r+1}$}}
\put(227.5,10){\makebox(0,0){$s_{n-2}$}}
\put(257.5,10){\makebox(0,0){$s_{n-1}$}}
\put(160,45){\makebox(0,0){$s_{n}$}}
\end{picture}
\caption[]{$T_{2,r,n-r}$ Dynkin diagram}
\end{figure}
 Together the reflections $s_1,\ldots,s_{n-1},s_n$
generate
the Weyl group with Dynkin diagram $T_{2,r,n-r}$, and the
inequality in Lemma \ref{cremona} is just what is needed to ensure
that this Weyl group is infinite.

It is clear that very general $n$-tuples of points in projective
space (that is, $n$-tuples outside of countably many proper subvarieties
of $(\P^{r-1})^n$)
are in Cremona general position. Thus Lemma \ref{cremona} gives
counterexamples to Hilbert's fourteenth problem over the complex
numbers (or any uncountable field). But this argument
leaves it unclear whether there
are any arrangements in Cremona general position
over the rational numbers or over finite fields. In the next
section we will construct such arrangements. These will occur in
the most interesting cases of Lemma \ref{cremona}, where
$\frac{1}{2}+\frac{1}{r}+\frac{1}{n-r}=1$, corresponding to
the affine Weyl groups of type $E_8^{(1)}$ ($n=9$, $r=3$ or 6)
or $E_7^{(1)}$ ($n=8$, $r=4$).

As a historical note, Coble knew in 1929 that the blow-up of $\P^2$ at
9 very general points contained infinitely many $(-1)$-curves
\cite[section 9]{Coble}, but the first proof to modern standards
was given by Nagata in 1960 \cite[Lemma 2.5 and Theorem 4a]{NagataRat2}.

\section{The method, with elementary proofs}
\label{elem}

In this section we prove Theorem \ref{main} over any sufficiently large 
finite field and over any infinite field.
This turns out to require only some simple geometry,
following the ideas of Coble and Dolgachev.

By Lemma \ref{cremona}, Theorem \ref{main} will follow if we can
produce explicit examples of $9$-tuples of points in $\P^2$ in
Cremona general position, or $8$-tuples in $\P^3$, or $9$-tuples
in $\P^5$. The following lemmas do that, over sufficiently large
finite fields and over infinite fields.
The idea is to use
arrangements of points which are special in one way (for example,
we take 9 points in $\P^2$ which are the intersection of two cubics),
but not too special.

In what follows, a {\it cubic }in $\P^2$ denotes an effective divisor
of degree 3; it is not assumed to be irreducible (likewise
for conics, quadrics in higher-dimensional projective spaces,
and so on). We say
that a cubic is {\it irreducible }if it is irreducible and reduced
(that is, the corresponding cubic form is irreducible).
A {\it pencil }of cubics in $\P^2$ means a linear system $\P^1$ of cubics;
equivalently, it is a 2-dimensional linear subspace of the vector
space of cubic forms in 3 variables. We say that two cubics
in $\P^2$ over a field $k$
intersect in 9 given points if those are all the intersection
points over the algebraic closure of $k$.

\begin{lemma}
\label{line}
Let $p_1,\ldots,p_9$ be 9 distinct rational points in $\P^2$ over a field $k$.
Suppose that there are cubics $C_1$ and $C_2$ whose intersection is the set
$\{p_1,\ldots, p_9\}$. Then no three of $p_1,\ldots,p_9$ lie on a line
if and only if all cubics in the pencil spanned by $C_1$ and $C_2$
are irreducible.
\end{lemma}

{\bf Proof. }Suppose that a cubic $C_1$ in the pencil is reducible.
Then we can write $C_1=L+D$ for some line $L$ and conic $D$. For any
other cubic $C_2$ in the pencil, 
$$\{p_1,\ldots,p_9\}=C_1\cap C_2 = (L\cap C_2)\cup (D\cap C_2).$$
Since $C_2$ intersects $C_1$ transversely,
three of these points are on the line $L$. Conversely,
suppose that $p_1,p_2,p_3$ lie on a line $L$. There is a unique cubic
on $L$ through $p_1,p_2,p_3$. Since we have a pencil (a $\P^1$)
of cubics in $\P^2$ through $p_1,p_2,p_3$, at least one cubic $C_1$
in the pencil must contain $L$. Thus $C_1$ is reducible. \qed

\begin{corollary}
\label{dim2}
Let $p_1,\ldots,p_9$ be 9 distinct rational points in $\P^2$ over a field $k$
which are the intersection of two cubics.
Suppose that no three of $p_1,\ldots,p_9$
lie on a line. Then $p_1,\ldots,p_9$ are in Cremona general position.
Therefore the blow-up $X$ 
of $\P^2$ at $p_1,\ldots,p_9$ has infinitely generated
total coordinate ring, and the corresponding 18-dimensional
representation of $(G_a)^6$ over $k$ has infinitely generated ring
of invariants.
\end{corollary}

The first part of Corollary \ref{dim2} was proved first by Manin
\cite[Lemma 3]{Manin}.

{\bf Proof. }The Cremona action on $\Pic (X)\cong \Z^{10}$ fixes
the anticanonical class $-K_X=3H-E_1-\cdots-E_9$. Therefore, performing
a standard Cremona transformation on any three  of $p_1,\ldots,p_9$
transforms a cubic through $p_1,\ldots,p_9$ into a cubic
through the new points $p_1',\ldots,p_9'$. Since $p_1,\ldots,p_9$
are the intersection of two cubics in $\P^2$, this remains true
after performing a standard Cremona transformation on any three of $p_1,\ldots,
p_9$. Since no three of $p_1,\ldots,p_9$ lie on a line, all cubics in the
pencil of cubics through these points are irreducible by Lemma \ref{line}.
The point is that this remains true of the proper transforms of these cubics on
the ``new'' $\P^2$. Therefore, no three of $p_1',\ldots,p_9'$ lie on a line.
Thus $p_1',\ldots,p_9'$ satisfy the same properties we assumed
for $p_1,\ldots,p_9$. We can repeat this process any number of times.
So $p_1,\ldots,p_9$ are in Cremona general position. Lemma \ref{cremona}
gives the rest. \qed

Corollary \ref{dim2} gives counterexamples to Hilbert's
fourteenth problem over all sufficiently large finite fields,
and over all infinite fields; we give explicit examples in
section \ref{allfields}.
The group involved is $(G_a)^6$. To get Mukai's best example,
with the group $(G_a)^3$, to work over the same fields, we can
use projective duality, as follows.

Let $q_1,\ldots,q_9$ be 9 distinct points on $\P^2$ which are the intersection
of two cubics. These points
can be represented by a linear map from $A^9$ onto $A^3$. The kernel
has dimension 6, and so dualizing
gives a linear map from $A^9$ onto a 6-dimensional vector space.
This gives the {\it dual }arrangement of 9 points $p_1,\ldots,p_9$ in
a 5-dimensional projective space.
Dolgachev gives an equivalent description of this arrangement,
under our assumption on $q_1,\ldots,q_9$:
$p_1,\ldots,p_9$ are the image of $q_1,\ldots,q_9$ under a Veronese
embedding $\P^2\arrow \P^5$ \cite[Prop.\ 5.4]{DolgachevEll}.

\begin{corollary}
\label{dim5}
Let $q_1,\ldots,q_9$ be 9 distinct rational points in $\P^2$ over a field $k$
which are the intersection of two cubics.
Suppose that no three of $q_1,\ldots,q_9$
lie on a line. Let $p_1,\ldots,p_9$ be the dual arrangement of
9 points in $\P^5$. Then $p_1,\ldots,p_9$ are in Cremona general position.
Therefore the blow-up $X$ 
of $\P^5$ at $p_1,\ldots,p_9$ has infinitely generated
total coordinate ring, and the corresponding 18-dimensional
representation of $(G_a)^3$ over $k$ has infinitely generated ring
of invariants.
\end{corollary}

{\bf Proof. }Since $q_1,\ldots,q_9$ are in linear general position
in $\P^2$, we check by hand that the dual arrangement $p_1,\ldots,p_9$ 
is in linear general position in $\P^5$. Performing the
standard Cremona transformation on $q_1,q_2,q_3$ in $\P^2$ and then dualizing
corresponds to performing the standard Cremona transformation
on $p_9,\ldots,p_4$ in $\P^5$, as Dolgachev-Ortland computed
\cite[Theorem VI.4]{DO}. (This is elementary: the assumption
that $q_1,\ldots,q_9$ are the intersection of two cubics
is not needed for this duality statement.)
Since $q_1,\ldots,q_9$ are in Cremona
general position in $\P^2$ by Corollary \ref{dim2}, it follows that
$p_1,\ldots,p_9$ are in Cremona general position in $\P^5$.
Lemma \ref{cremona} gives the rest. \qed

Thus we have counterexamples to Hilbert's fourteenth problem
using an 18-dimensional representation of $(G_a)^3$ over all
sufficiently large finite fields and all infinite fields.
Finally, let us check that Mukai's
final class of examples, using the group $(G_a)^4$, also works
over all sufficiently large finite fields and all infinite fields.
A {\it net }of quadrics in $\P^3$ means a linear system
$\P^2$ of quadrics, or equivalently a 3-dimensional linear subspace
of the vector space of quadratic forms in 4 variables.

\begin{lemma}
\label{plane}
Let $p_1,\ldots,p_8$ be 8 distinct points in $\P^3$ which are
the intersection of three quadrics, $Q_1\cap Q_2\cap Q_3$.
Then no four of the points $p_1,\ldots,p_8$
lie on a plane if and only if every quadric in the net spanned
by $Q_1,Q_2,Q_3$ is irreducible.
\end{lemma}

{\bf Proof. }If some quadric $Q$ in the net is reducible, then
$Q=S_1\cup S_2$ for some planes $S_1$ and $S_2$ in $\P^3$. Then four of the
points $p_1,\ldots,p_8$ lie on $S_1$ (and the other four lie on $S_2$).

Conversely, suppose that $p_1,p_2,p_3,p_4$ lie on a plane $S$ in $\P^3$.
No three of the points $p_1,\ldots,p_8$ lie on a line $L$ in $\P^3$;
otherwise every quadric in the net would contain three points 
on $L$, hence would contain $L$, contradicting that
$Q_1\cap Q_2\cap Q_3=\{p_1,\ldots,p_8\}$. So no three of the points
$p_1,p_2,p_3,p_4$ lie on a line. Therefore
these four points are the complete intersection of two conics in the plane $S$,
$p_1p_2\cup p_3p_4$ and $p_1p_3\cup p_2p_4$. It follows that there
is only a pencil (a $\P^1$) of conics in $S$ through $p_1,p_2,p_3,p_4$.
Since we have a net (a $\P^2$) of quadrics in $\P^3$ that contain
$p_1,p_2,p_3,p_4$, at least one quadric $Q$ in the net must contain
the whole plane $S$. Thus $Q$ is reducible. \qed

\begin{corollary}
\label{dim3}
Let $p_1,\ldots,p_8$ be 8 distinct points in $\P^3$ which are
the intersection of three quadrics, and suppose that
no four of the points $p_1,\ldots,p_8$
lie on a plane. Then $p_1,\ldots,p_8$ are in Cremona general
position. Therefore the blow-up $X$ 
of $\P^3$ at $p_1,\ldots,p_8$ has infinitely generated
total coordinate ring, and the corresponding 16-dimensional
representation of $(G_a)^4$ over $k$ has infinitely generated ring
of invariants.
\end{corollary}

The first part of Corollary \ref{dim3} was apparently known
to Coble \cite[section 44, last paragraph]{Coble}.

{\bf Proof. }The Cremona action on $\Pic (X)$ fixes half the anticanonical
class, $-\half K_X=2H-E_1-\cdots-E_8$. So performing a standard
Cremona transformation on any four of $p_1,\ldots,p_8$ transforms
a quadric through $p_1,\ldots,p_8$ into a quadric through the new
points $p_1',\ldots,p_8'$. Since $p_1,\ldots,p_8$ are the intersection
of three quadrics, $p_1',\ldots,p_8'$ are also the intersection
of three quadrics. By Lemma \ref{plane}, since no four of $p_1,
\ldots,p_8$ lie on a plane, all the quadrics through $p_1,\ldots,p_8$
are irreducible. Therefore their proper transforms, the quadrics
through $p_1',\ldots,p_8'$, are also all irreducible. So no four
of $p_1',\ldots,p_8'$ lie on a plane. Thus $p_1',\ldots,p_8'$
satisfy the same assumptions as $p_1,\ldots,p_8$, and so we can
repeat the process any number of times. That is, $p_1,\ldots,p_8$
are in Cremona general position. \qed

\section{Examples over arbitrary fields}
\label{allfields}

We now prove Theorem \ref{main}, showing that our three classes
of representations with infinitely generated rings of invariants
all exist over arbitrary fields. In later sections, we will
describe the geometry behind these representations in more detail.
One benefit will be to construct representations 
with infinitely generated ring of invariants that are
given by simpler formulas. Another benefit will be to prove
partial results about exactly
when finite generation holds.
For now, we will just prove
Theorem \ref{main} as stated.

Theorem \ref{mukai} associates a $2n$-dimensional representation
of $(G_a)^{n-r}$ over a field $k$ to any arrangement of $n$ distinct
$k$-points of $\P^{r-1}$ not contained in a hyperplane. More generally,
we can define such a representation associated to a smooth 0-dimensional
subscheme of degree $n$ defined over $k$, even if the individual
points are not defined over $k$, using the standard technique
of ``twisting'' in Galois cohomology \cite[III.1.3]{SerreGal}.
Explicitly, given a Galois extension
field $K/k$ over which the points are defined, let the Galois group
$\Gal(K/k)$ act on $K^n$ by permuting the basis vectors as $\Gal(K/k)$
permutes the $n$ points, and by $\sigma(ax)=\sigma(a)\sigma(x)$ for
$a\in K$. We can view $K^n$ with this Galois action as a ``twisted
form'' $H$ of the group $(G_a)^n$ over $k$, but it is in fact isomorphic
to $(G_a)^n$ over $k$ because $H(k)=(K^n)^{\Gal(K/k)}$ is an $n$-dimensional
vector space over $k$. A linear map $K^n\arrow K^r$ associated
to the $n$ points in $\P^{r-1}$ is evidently Galois-equivariant, and
so (taking Galois invariants) it gives a linear map $H(k)\cong k^n\arrow
k^r$. The kernel is a vector space over $k$, which canonically
determines an algebraic group $G$ over $k$, clearly
isomorphic to $(G_a)^{n-r}$.  Likewise, the $2n$-dimensional
representation of $(G_a)^{n}$ in Theorem \ref{mukai} 
gives a $2n$-dimensional representation of its twisted form
$H$ and hence of the subgroup $G\cong (G_a)^{n-r}$ over $k$.

Moreover, whether a ring of invariants $O(V)^G$ for a representation $V$ of
an algebraic group $G$ is finitely generated does not change under
extension of the base field, since the ring of invariants over
an extension field $K/k$ is just the ring of invariants over $k$
tensored with $K$.
So we can apply the previous section's results. For example, given
two plane cubics $C_1$ and $C_2$
over a field $k$ whose intersection is smooth of dimension 0,
consider the associated 18-dimensional representation of
$(G_a)^6$ over $k$.
If no three points of $C_1\cap C_2$ over the algebraic closure
of $k$ lie on a line, then the corresponding representation of 
$(G_a)^6$ over $k$ has infinitely generated ring of invariants,
by Corollary \ref{dim2}.
In practice, it is often easier to check the equivalent hypothesis that
all cubics in the pencil spanned by $C_1$ and $C_2$ are irreducible.

For example, consider the cuspidal cubics $(y+z)^3+xz^2$ and $x^3+y^2z$
over a field $\F_p$. For primes $p\leq 23$ (not an optimal bound), 
we compute that the intersection of these two
cubics in $\P^2$ is smooth of dimension 0, and that all nonzero
linear combinations
of these two cubics are irreducible. (This is particularly easy to check
for $p=2$ or 3, where the given cuspidal cubics are the
only singular cubics in the pencil.) By Corollary \ref{dim2},
for $p\leq 23$, the
associated 18-dimensional representation of $(G_a)^6$ over $\F_p$
has infinitely generated ring of invariants.
For fields of characteristic greater
than 23 (or characteristic zero),
we can apply Corollary \ref{dim3} directly. Namely, the
9 points in $\P^2$ which form the columns of the following matrix
are the base locus of a pencil of cubics, and are in linear
general position over $\Q$ and over $\F_p$ for all $p>23$. (One can look
for such examples by choosing eight points
in $\P^2$ and computing
the ninth point on the pencil of cubics through the eight points.)
$$\begin{pmatrix}
            1 & 0 & 0 & 1 &  1 &  2 & -3 & -2 & -7\\
            0 & 1 & 0 & 1 & -1 & -1 &  4 & -5 &  2\\
            0 & 0 & 1 & 1 &  2 &  1 &  1 &  1 & -1\\
\end{pmatrix}$$
Thus we have constructed 18-dimensional
representations of $(G_a)^6$ with infinitely generated ring of invariants
over an arbitrary field.

Now consider again the intersection of the 
cuspidal cubics $(y+z)^3+xz^2$ and $x^3+y^2z$
over a field $\F_p$ with $p\leq 23$. The projective dual of these
9 points in $\P^2$ gives a smooth subscheme of degree 9 in $\P^5$
over $\F_p$. Consider the associated
representation of $(G_a)^3$ over $\F_p$.
Since all cubics in this pencil of cubics in $\P^2$
are irreducible when $p\leq 23$,
Corollary \ref{dim5} shows that the ring of invariants for
this 18-dimensional
representation of $(G_a)^3$ over $\F_p$ is infinitely generated.
For fields of characteristic greater than $23$ or of characteristic zero,
we can apply Corollary \ref{dim5}
directly, using the previous paragraph's arrangement of 9 rational points.
Thus we have constructed 18-dimensional
representations of $(G_a)^3$ with infinitely generated ring of invariants
over an arbitrary field.

We now construct analogous representations of $(G_a)^4$.
Over a field $\F_p$, consider
the quadrics $xy+y^2+z^2$, $xw+y^2+w^2$, $xz-zw+z^2+w^2$.
We compute that for all primes $p\leq 7$ (not an optimal bound),
the intersection of these three quadrics is smooth of dimension 0,
and all nonzero linear combinations of these quadrics are irreducible.
Therefore, for $p\leq 7$, the corresponding
16-dimensional representation of $(G_a)^4$ over $\F_p$ has infinitely
generated ring of invariants. For fields of characteristic greater
than 7 (or characteristic zero),
we can apply Corollary \ref{dim3} directly. Namely, the
8 points in $\P^3$ which form the columns of the following matrix
are the base locus of a net of quadrics, and are in linear
general position over $\Q$ and over $\F_p$ for all $p>7$. (One can look
for such examples by choosing seven points $p_1,\ldots,p_7$
in $\P^3$ and computing
the eighth point on the net of quadrics through $p_1,\ldots,p_7$.)
$$\begin{pmatrix}
            1 & 0 & 0 & 0 & 1 & -2 &  1  & -6\\
            0 & 1 & 0 & 0 & 1 &  2 & -4 & -8\\
            0 & 0 & 1 & 0 & 1 &  1 & -3 & -7\\
            0 & 0 & 0 & 1 & 1 & -3 &  4 & -4
\end{pmatrix}$$
Thus we have proved the existence of 16-dimensional
representations of $(G_a)^4$ with infinitely generated ring of invariants
over an arbitrary field.
 \qed (Theorem \ref{main})

\section{Elliptic fibrations, and representations of $(G_a)^6$ on $A^{18}$}

We now describe the geometry behind the representations of
$(G_a)^6$ in Theorem \ref{main}. It turns out that finite generation
depends upon the Mordell-Weil group of a certain
elliptic fibration of $\P^2$.
For a certain class of representations, we can say exactly
when the ring of invariants is finitely generated and when
it is not. As a concrete application, we give examples
of non-finite generation where the coefficients
of the representation are very simple (in particular,
simpler than the examples in section \ref{allfields}). 
As in section \ref{allfields}, we get examples of non-finite
generation over all fields,
even the field of order 2.

The general problem here is to understand
the border between finite generation and non-finite generation
(for rings of invariants, or for total coordinate rings). Section
\ref{elem} shows that if we assume a small amount of general position, we get
examples of non-finite generation. In this section we find that
even for some more special arrangements of points, which do not satisfy
the hypotheses of section \ref{elem}, we can still prove
non-finite generation. This case of 9 points in the plane has been
intensely studied, and many of the results of this section
can be deduced from various earlier works, as we will see.
The main novelty is the precise characterization
of finite generation in Theorem \ref{sharp2}. 

\begin{corollary}
\label{six}
Let $(G_a)^n$ act on $V=A^{2n}$ by
$$(t_1,\ldots,t_n)(x_1,\ldots,x_n,y_1,\ldots,y_n)=(x_1,\ldots,x_n,
y_1+t_1x_1,\ldots,y_n+t_nx_n).$$
Over any field $k$ of characteristic not 2 or 3,
consider the subgroup $G=(G_a)^6$ of $(G_a)^9$ which is the kernel
of the following linear map $A^9\arrow A^3$:
$$\begin{pmatrix}
           -1 & -1 & -1 & 0 & 0 & 0 & 1 & 1 & 1\\
           -1 &  0 & 1  & -1& 0 & 1 & -1& 0 & 1\\
            1 &  1 & 1  & 1 & 1 & 1 &  1& 1 & 1
\end{pmatrix}$$
Consider the restriction of the above 18-dimensional representation
of $(G_a)^9$ to the subgroup $(G_a)^6$.
The ring of invariants of this representation
is not finitely generated over $k$.

There are other (almost equally simple) 18-dimensional
representations of $(G_a)^6$ over $\F_2$ or $\F_3$ for which
the ring of invariants is not finitely generated.
\end{corollary}

Corollary \ref{six}
is a consequence of the following theorem, which characterizes
exactly which intersections of two cubics in $\P^2$
yield infinitely generated rings of invariants (or infinitely generated
total coordinate rings). We generalize this theorem
to all rational elliptic surfaces in Theorem \ref{rational}.

\begin{theorem}
\label{sharp2}
Let $p_1,\ldots,p_9$ be 9 distinct rational points in $\P^2$ over a field $k$
which are the intersection of two cubics.
Let $a$ be the number of collinear triples
of points $p_1,\ldots,p_9$, let $b$ be
the number of partitions of $p_1,\ldots,p_9$
into three collinear triples, and let
$$\rho=8-a+b.$$
Then $\rho\geq 0$.

If $\rho>0$, then
the blow-up $X$ 
of $\P^2$ at $p_1,\ldots,p_9$ has infinitely generated
total coordinate ring, and the corresponding 18-dimensional
representation of $(G_a)^6$ over $k$ has infinitely generated ring
of invariants.

Conversely, if $\rho=0$, then these rings
are finitely generated.
\end{theorem}

Theorem \ref{sharp2} strengthens Corollary \ref{dim2},
which proves infinite generation assuming that there are no collinear
triples among $p_1,\ldots,p_9$.

{\bf Proof. }Let $X$ be the blow-up of $\P^2$ at the points $p_1,\ldots,p_9$.
The pencil of cubics through $p_1,\ldots,p_9$ makes 
the blow-up $X$ an elliptic surface $f:X\arrow \P^1$
(or perhaps quasi-elliptic, in characteristics 2 and 3).
We want to give criteria for $X$ to contain infinitely many $(-1)$-curves,
without requiring that no three of $p_1,\ldots,p_9$ lie on a line
(as in Corollary \ref{dim2}). There is in fact a complete classification
of the intersections of two plane cubics such that
the blow-up $X$ contains only finitely many $(-1)$-curves ($X$
is called an {\it extremal }rational elliptic surface),
by Miranda-Persson in characteristic zero and W.~Lang in positive
characteristic \cite{MP, Lang1, Lang2}; see also Cossec-Dolgachev
\cite[section 5.6]{CD}. There is also more general work by Nikulin
classifying certain types of varieties whose cone of curves
is finite polyhedral \cite{Nikulin1, Nikulin2}.
We will not use any of
these classification results, but just explain how to check
whether a given rational elliptic surface $X$ has infinitely many
$(-1)$-curves.

The $(-1)$-curves in any surface
can be described as the smooth rational curves
$C$ such that $(-K_X)\cdot C=1$. Since the fibers of the elliptic
fibration $f:X\arrow \P^1$ are in the linear system $|-K_X|$, it follows
that the $(-1)$-curves in $X$ are precisely the sections of the
elliptic fibration. Thus, to make $X$ have infinitely many 
$(-1)$-curves, it suffices to arrange that the Mordell-Weil
group $\Pic^0(E)$ of the general fiber (an elliptic curve over the field
$k(t)$) has rank at least 1. We can define the Mordell-Weil
group as the group of sections of this elliptic fibration,
with one section considered as the zero section.

For each reducible fiber $F$ of the elliptic fibration $X\arrow \P^1$,
let $r_F+1$ be the number of irreducible components of $F$. 
By Tate \cite[4.5]{Tate} and 
Shioda \cite[Cor.\ 1.5]{Shioda}, the Mordell-Weil rank of $E$
over $k(t)$ is equal to $8-\sum r_F$. For the reader's convenience,
here is a proof. The Picard group of the general fiber $E$ over $k(t)$
is the quotient of the Picard group of the surface $X$ by
the classes of all irreducible divisors in $X$ which do not map onto $\P^1$.
So $\Pic(E)$ is the quotient
of $\Pic(X)\cong \Z^{10}$ by the class of $-K_X$ (the class of each irreducible
fiber of $f:X\arrow \P^1$) and by
$r_F$ classes for each fiber $F$ with $r_F+1$ irreducible
components. Moreover, if we write the reducible fibers as
$F=\sum_{j=1}^{r_F+1}m_{Fj}D_{Fj}$, then the divisors $D_{Fj}$
for all $F$ and all $1\leq j\leq r_F$, together with $-K_X$,
are linearly independent in $\Pic(X)_{\Q}$ (even modulo numerical equivalence).
This is a fact about any morphism from a surface to a curve
\cite[Cor.~VIII.4]{Beauville}. Therefore $\Pic(E)_{\Q}$ has rank
$10-1-\sum r_F$ and the Mordell-Weil group $\Pic^0(E)_{\Q}$
has rank $8-\sum r_F$, as we want.

The only possible reducible cubics in the given pencil are
the sum of a line and a conic, or the sum of three lines. By the proof
of Lemma \ref{line}, the lines occurring in these reducible
cubics are precisely the lines through collinear triples among
the points $p_1,\ldots,p_9$.
Therefore $\sum r_F$ is equal to the number $a$ of collinear triples,
minus the number $b$
of partitions of $p_1,\ldots,p_9$ into three disjoint collinear
triples. (The second term comes in because a cubic in the pencil
which consists of three lines contributes only 2, not 3, to
$\sum r_F$.) So the Mordell-Weil rank of $E$ over $k(t)$
is $\rho=8-a+b$, which is therefore nonnegative. If $\rho$ is greater
than 0,
then the blow-up $X$ contains infinitely many $(-1)$-curves. It follows
that the total coordinate ring of $X$ is infinitely generated, and that
the corresponding ring of invariants for an 18-dimensional
representation of $(G_a)^6$ over $k$ is infinitely generated.

It remains to show that when $\rho=0$, then the total coordinate
ring of $X$ is finitely generated. Our main tools will be the
cone and contraction theorems, proved for log surfaces in any characteristic
by Koll\'ar and Kov\'acs \cite[2.1.1, 2.3.3]{KK},
and Hu and Keel's notion of a Mori dream space \cite{HK}. In particular,
for any klt pair $(X,\Delta)$ of dimension 2, we can contract any
$(K_X+\Delta)$-negative extremal face of the closed cone
of curves $\overline{NE}(X)$, yielding a projective
variety. Koll\'ar
and Kov\'acs only state this for contractions of extremal rays,
but the result for higher-dimensional extremal faces is a consequence.
The point is that for
any $(K_X+\Delta)$-negative extremal face $F$, if the
contraction $f:X\arrow Y$ of one ray of $F$
is birational, then the pair $(Y,f_*\Delta)$ is klt and
the image of $F$ is a
$(K_Y+f_*\Delta)$-negative extremal face \cite[Lemma 2.3.5]{KK}.
Then the cone theorem (including the contraction theorem)
for higher-dimensional extremal faces
follows by induction on the Picard number of $X$.

The assumption $\rho=0$ tells us
that the subspace $(K_X)^{\perp}\cong \Q^9$ of $\Pic(X)\otimes\Q$ is spanned
by the curves in $X$ which are contained in fibers of the morphism
$X\arrow \P^1$. The sum of all the curves in a reducible
fiber, with multiplicities, is numerically equivalent to a general fiber
$-K_X$; so $(K_X)^{\perp}$ is spanned just by the finitely many
irreducible components of reducible fibers of $X\arrow \P^1$.
We can also describe these curves as the $(-2)$-curves on $X$
(curves isomorphic to $\P^1$ with self-intersection $-2$), using the
fact that $-K_X$ has degree 0 on any $(-2)$-curve.

Another consequence of the assumption $\rho=0$ is that the Mordell-Weil
rank of the elliptic surface $X\arrow \P^1$ is zero. Since the sections of
$X\arrow \P^1$ are exactly the $(-1)$-curves, this means that $X$
contains only finitely many $(-1)$-curves.  By the cone theorem,
the extremal rays of the closed cone of curves
in the $K_X$-negative half space $(K_X)^{<0}$ are all spanned
by $(-1)$-curves, and so there are only finitely many
extremal rays in $(K_X)^{<0}$. Since $-K_X$ is semi-ample
(corresponding to the contraction $X\arrow \P^1$), we know
that the closed cone of curves is contained in $(K_X)^{\leq 0}$.

Moreover, for any curve $C$ on $X$ which is not one of the
finitely many $(-2)$-curves $D_{Fj}$, we have $C\cdot D_{Fj}\geq 0$
for all $F$ and $j$. Because the curves $D_{Fj}$ span $(K_X)^{\perp}$
and are divided into subsets with $\sum_j m_{Fj}D_{Fj}\sim -K_X$,
where $m_{Fj}>0$, the cone 
$$\{ x\in N^1(X):x\cdot D_{Fj}\geq 0\text{ for all }F,j\}$$
is a finite polyhedral subcone of $(K_X)^{\leq 0}$ whose intersection
with $(K_X)^{\perp}$ is just $\R^{\geq 0}\cdot (-K_X)$. Since all curves except
the $(-2)$-curves belong to this cone, it follows that the closed
cone of curves is finite polyhedral, spanned by the
$(-1)$-curves and $(-2)$-curves. (This consequence of
$\rho=0$ was proved already by Nikulin \cite[Example 1.4.1]{Nikulin1}.)

By Hu and Keel, a projective variety has finitely generated
total coordinate ring if and only if it is
a ``Mori dream space'' \cite[Proposition 2.9]{HK}.
For a smooth projective surface, this means that
the first Betti number is zero, the closed
cone of curves is finite polyhedral, and every codimension-1
face of the cone can be contracted.

Thus, it remains to show that every codimension-1 face of the closed
cone of curves $\overline{NE}(X)$
can be contracted. A face in the $K_X$-negative
half space (equivalently, a face spanned by $(-1)$-curves) can be contracted
by the cone theorem. On the other hand, because the $(-2)$-curves
span the hyperplane $(-K_X)^{\perp}$, the
intersection $\overline{NE}(X)\cap (K_X)^{\perp}$ is a codimension-1
face (the span of all the $(-2)$-curves). We can contract
this face using the line bundle $-K_X$, corresponding
to the contraction $X\arrow \P^1$.

It remains to contract a codimension-1 face $A$ of $\overline{NE}(X)$
spanned by some $(-1)$-curves and some $(-2)$-curves. Let
$C_1,\ldots,C_r$ be the $(-2)$-curves in the face $A$. These
cannot include all the irreducible components of any fiber $F$
of $X\arrow \P^1$, because then $-K_X\sim \sum_j m_{Fj}D_{Fj}$
would belong to the face $A$. But $-K_X$ is a positive linear
combination of all the $(-2)$-curves. Therefore $-K_X$
belongs to the interior
of the codimension-1 face $(K_X)^{\perp}\cap \overline{NE}(X)$,
and so it cannot belong to any other face.

Since the $(-2)$-curves $C_1,\ldots,C_r$ do not include
all the irreducible components of any fiber of $X\arrow \P^1$,
the intersection pairing on $\oplus_{i=1}^r \Q C_i$
is negative definite \cite[Cor.~VIII.4]{Beauville}. So there
are rational numbers $a_i$ such that $(\sum a_iC_i)\cdot C_j=-1$
for all $j=1,\ldots,r$. By Koll\'ar-Mori \cite[Lemma 3.41]{KM},
it follows that $a_i\geq 0$ for all $i$.

Let
$\Delta=\epsilon\sum a_iC_i$ for some small positive number
$\epsilon$. The pair $(X,\Delta)$ is klt \cite[Definition 2.2.6]{KK}
% \cite[Definition 2.34]{KM}
since $X$ is smooth and $\epsilon$ is small; so we can apply
the contraction theorem. Clearly $K_X+\Delta$ is negative on the
$(-1)$-curves $C$ in the face $A$, since $K_X\cdot C=-1$
and $\epsilon$ is small. Also, $K_X+\Delta$ has degree $-\epsilon $
on the $(-2)$-curves $C_1,\ldots,C_r$ in $A$. Thus the face $A$
is $(K_X+\Delta)$-negative. By the cone theorem,
we can contract the face $A$. This completes the proof
that the total coordinate
ring of $X$ is finitely generated.
\qed (Theorem \ref{sharp2})

We now prove Corollary \ref{six}.
For any field $k$ of characteristic not 2 or 3,
consider the 9 points in $\P^2(k)$ given in affine coordinates by
$(x,y)$ where $x$ and $y$ run through the set $\{-1,0,1\}$.
These points are the intersection of the two cubics $x^3=xz^2$
and $y^3=yz^2$, each consisting of three lines through a point in $\P^2$.
We count (using the assumption on the characteristic) that
there are 8 collinear triples among these points $p_1,\ldots,p_9$,
and 2 partitions of $p_1,\ldots,p_9$ into three collinear triples.
(The two partitions into three collinear triples
correspond to two cubics $x^3-xz^2=x(x-z)(x+z)$ and $y^3-yz^2=y(y-z)(y+z)$
in our pencil
which are unions of three lines, and the other two collinear triples
correspond to two other reducible cubics in the pencil,
$(x^3-xz^2)-(y^3-yz^2)=(x-y)(x^2+xy+y^2-z^2)$ and
$(x^3-xz^2)+(y^3-yz^2)=(x+y)(x^2-xy+y^2-z^2)$.)
So $\rho=8-8+2=2$ is greater than 0. By Theorem \ref{sharp2},
the blow-up $X$ of $\P^2$ at this
set of 9 points has infinitely many $(-1)$-curves. Thus we get
an 18-dimensional representation of $(G_a)^6$ over
$\F_p$ for $p\geq 5$ (or over $\Q$) whose ring of invariants
is not finitely generated. More precisely, the Mordell-Weil group
has rank 2 in this case.

We refer to the Appendix for simple examples of 18-dimensional
representations of $(G_a)^6$ over $\F_2$ or $\F_3$
whose rings of invariants are not finitely generated.
\qed (Corollary \ref{six})

\section{A new fibration of $\P^5$ by abelian surfaces, and
representations of $(G_a)^3$ on $A^{18}$}
\label{abelian}

We now give a richer geometric explanation for the infinitely
generated rings of invariants for $(G_a)^3$ constructed in
Corollary \ref{dim5}: they are explained by a fibration of
a blow-up of $\P^5$ by abelian surfaces.

\begin{theorem}
\label{three}
Let $q_1,\ldots,q_9$ be nine distinct
rational points in $\P^2$ over a field $k$
which are the intersection of two cubics. Suppose that no three
of $p_1,\ldots,p_9$ lie on a line. Let $p_1,\ldots,p_9$ be
the dual arrangement of 9 points in $\P^5$. Then the blow-up
$X$ of $\P^5$ at $p_1,\ldots,p_9$ is pseudo-isomorphic to
a smooth projective variety $W$ which is an abelian surface
fibration over $\P^3$ with a section.
The Mordell-Weil group of $W$ over $\P^3$ has rank 8. The translates
by the Mordell-Weil group of the $(-1)$-divisors $E_1,\ldots,E_9$
yield infinitely many $(-1)$-divisors on $W$, or equivalently on $X$.
As a result,
the total coordinate ring of $X$ is not finitely generated,
and the corresponding 18-dimensional representation of $(G_a)^3$
over $k$ has infinitely generated ring of invariants.
\end{theorem}

{\bf Proof. }We recall Dolgachev's equivalent description
of the dual arrangement,
under our assumption on $q_1,\ldots,q_9$:
$p_1,\ldots,p_9$ are the image of $q_1,\ldots,q_9$ under a Veronese
embedding $\P^2\arrow \P^5$ \cite[Prop.\ 5.4]{DolgachevEll}.

There is a natural homomorphism from the group $\Z^8\rtimes \Z/2$
to the automorphism group of the blow-up $Y$ of $\P^2$ at
$q_1,\ldots,q_9$ \cite[p.~124]{DO}. Explicitly, $Y$ is an
elliptic
surface, with the $\Z^8$ subgroup giving translations by differences
of sections, and the $\Z/2$ giving the map $z\mapsto -z$ with respect
to some zero-section. The construction of $p_1,\ldots,p_9$ by duality
gives a corresponding action of $\Z^8\rtimes \Z/2$ by
pseudo-automorphisms of the blow-up
$X$ of $\P^5$ at $p_1,\ldots,p_9$, by Dolgachev \cite[5.12]{DolgachevEll}.
Moreover,
these pseudo-automorphisms are defined on the blown-up Veronese surface
$Y$ inside $X$ and give the action we mentioned of $\Z^8\rtimes \Z/2$
on $Y$. Our goal is to give a more geometric
interpretation of these pseudo-automorphisms of $X$.

The linear system of
$-\half K_X = 3H-2E_1-\cdots -2E_9$ (that is, cubics on $\P^5$
singular at $p_1,\ldots,p_9$) is a $\P^3$, with base locus
a union of 45 curves: the 36 lines through pairs of points $p_i$
and the 9 rational normal curves through any 8 of the points
$p_1,\ldots,p_9$. Indeed, the base locus of the linear system
$|-\half K_X|$, viewed on $\P^5$,
is an intersection
of 4 cubics. One checks
that the base locus has dimension 1, which is easy in a particular
example (and therefore holds for a general 9-tuple $p_1,\ldots,p_9$
as above). Thus the base locus is a complete intersection
of four cubics in $\P^5$ and hence has degree 81.
The 45 curves mentioned have total degree 81,
and they are all in the base locus, because
$-\half K_X$ has degree $-1$ on these curves. So
the base locus is precisely the union of these 45 curves,
each with multiplicity 1.

We can describe the rational map $X\dashrightarrow \P^3$
associated to $-\half K_X$. The base locus
of this linear system consists of 45 disjoint curves isomorphic to $\P^1$
in $X$,
all with normal bundle $O(-1)^{\oplus 4}$. Using the linear
system $|-\half K_X|$, we can perform an inverse flip on these
45 curves, giving another smooth projective 5-fold $W$
in which the 45 $\P^1$'s with normal bundle $O(-1)^{\oplus 4}$
have been replaced by $\P^3$'s with normal bundle $O(-1)^{\oplus 2}$.
Explicitly, $W$ can be described as $\text{Proj}$ of the 
ring $R(X,L)=\oplus_{a\geq 0} H^0(X,aL)$ for any line bundle $L$ which
is a positive linear combination of $H+5(-\half K_X)$
and $-\half K_X$.
On this pseudo-isomorphic variety $W$, the linear system
$|-\half K_W|$ ($=|-\half K_X|$) is basepoint-free, giving
a morphism $W\arrow \P^3$. The generic fiber of $W\arrow \P^3$
is a principally polarized abelian surface. The
$(-1)$-divisors $E_1,\ldots,E_9$ on $W$ all induce the same polarization
of the generic fiber, which is two times a principal
polarization. (Geometrically, each divisor $E_i$ intersects
the general fiber in a curve of genus 5.) The 45 $\P^3$'s we
have produced in $W$ are sections of the fibration, and their differences
generate only
a rank 8 subgroup of the Mordell-Weil group,
for $p_1,\ldots,p_9$ general as above.

To check that the smooth fibers of $W\arrow \P^3$ are abelian surfaces,
note that the canonical bundle $K_W$ is trivial on all fibers,
since $-\half K_W$ is pulled back from $\P^3$. So the smooth fibers
have trivial canonical bundle and hence (say, in characteristic not 2 or 3)
are either K3 surfaces or abelian surfaces. To see that they are
abelian surfaces, it suffices to compute that $c_2(W)(-\half K_W)^3=0$,
since K3 surfaces have Euler characteristic $c_2$ equal to 24
while abelian surfaces have $c_2=0$. To do that calculation,
one can compute that $c_2(X)(-\half K_X)^3=45$ and that performing
an inverse flip (replacing a $\P^1$ with normal bundle
$O(-1)^{\oplus 4}$ by a $\P^3$ with normal bundle $O(-1)^{\oplus 2}$)
lowers $c_2(X)(-\half K_X)^3$ by 1.

The following result helps to show what is going on.

\begin{lemma}
\label{minimal}
Let $f:X\arrow S$ be a proper morphism of smooth varieties
over a field such that
$K_X$ has degree zero on all curves $C$
with $f(C)=$ point. Suppose that the smooth locus of the generic
fiber has a group structure (for example, the generic fiber
could be an abelian variety, or a cuspidal cubic curve in characteristic
2 or 3). Then
for any sections $S_1$ and $S_2$ of $f$ over the generic point
of $S$, there is a pseudo-automorphism of $X$ over $S$
that maps $S_1$ to $S_2$.
\end{lemma}

{\bf Proof. }We use a basic fact of minimal model theory:
for any proper morphisms $f_i:X_i\arrow S$, $i=1,2$,
such that $X_i$ is smooth and $K_{X_i}$ has nonnegative
degree on all curves that map to a point in $S$, every birational map
from $X_1$ to $X_2$ over $S$ is a pseudo-isomorphism. The proof
works in any characteristic \cite[Theorem 3.52]{KM}, and applies more generally
to varieties with terminal singularities. Any section of the
morphism $f:X\arrow S$ over the generic point of $S$ automatically lies
in the smooth locus of the generic fiber, which we assume has
a group structure. Therefore,
for any sections $S_1$ and $S_2$ over the generic point,
adding $S_2-S_1$ using the group structure
is a pseudo-automorphism of $X$ over $S$. \qed

Since the Mordell-Weil group of the generic fiber of $W\arrow 
\P^3$ has rank 8,
it is in particular infinite. By the fact of minimal model theory
just stated, the Mordell-Weil group acts on $W$, or equivalently
on $X$, by pseudo-automorphisms. These pseudo-automorphisms
of $X$ map the $(-1)$-divisor $E_1$ to infinitely many other
$(-1)$-divisors. It follows that the total coordinate ring of $X$
is not finitely generated, and that the corresponding representation
of $(G_a)^3$ has infinitely generated ring of invariants, as we already
knew from Corollary \ref{dim5}.
\qed (Theorem \ref{three}).

\section{An elliptic fibration of $\P^3$, and
representations of $(G_a)^4$ on $A^{16}$}

We now describe the geometry behind the representations of
$(G_a)^4$ in Theorem \ref{main} in more detail.
As in Theorem \ref{sharp2},
the key concept is the Mordell-Weil group of an elliptic fibration,
in this case an elliptic fibration of a blow-up of $\P^3$.
For a certain class of representations, we can show that
the ring of invariants is infinitely generated except in
a very special situation. As a concrete application, we give examples
of non-finite generation where the coefficients
of the representation are very simple (in particular,
simpler than the examples in section \ref{allfields}). 
As in section \ref{allfields}, we get examples of non-finite
generation over all fields,
even the field of order 2.

\begin{corollary}
\label{four}
Over any field $k$ of characteristic not 2,
start with Nagata's representation of $(G_a)^8$ on
$A^{16}$ as in Theorem \ref{mukai}, and restrict to the subgroup $(G_a)^4$
spanned by the rows of the following $4\times 8$ matrix.
$$\begin{pmatrix}
            0 & 0 & 0 & 0 & 1 & 1 & 1 & 1\\
            0 & 0 & 1 & 1 & 0 & 0 & 1 & 1\\
            0 & 1 & 0 & 1 & 0 & 1 & 0 & 1\\
            1 & 1 & 1 & 1 & 1 & 1 & 1 & 1
\end{pmatrix}$$
Then the ring of invariants is not finitely generated over $k$.
There are other (almost equally simple) 
16-dimensional
representations of $(G_a)^4$ over $\F_2$ for which
the ring of invariants is not finitely generated.
\end{corollary}

Corollary \ref{four} is a consequence
of the following theorem. Combined with A.~Prendergast-Smith's
results discussed below, Theorem \ref{sharp3} characterizes
exactly which intersections of three quadrics in $\P^3$
yield infinitely generated rings of invariants (or infinitely generated
total coordinate rings). 

\begin{theorem}
\label{sharp3}
Let $p_1,\ldots,p_8$ be 8 distinct rational points in $\P^3$ over a field $k$
which are the intersection of three quadrics.
Let $a$  be the number of coplanar quadruples
of points $p_1,\ldots,p_8$, and define
$$\rho=7-\frac{a}{2}.$$
Then $\rho\geq 0$.
If $\rho$ is greater than 0, then 
the blow-up $X$ 
of $\P^3$ at $p_1,\ldots,p_8$ has infinitely generated
total coordinate ring, and the corresponding 16-dimensional
representation of $(G_a)^4$ over $k$ has infinitely generated ring
of invariants.
\end{theorem}

Theorem \ref{sharp3} strengthens Corollary \ref{dim3},
which proves infinite generation assuming that there are no coplanar
quadruples among $p_1,\ldots,p_8$.
A.~Prendergast-Smith has classified the intersections
of three quadrics in $\P^3$ with $\rho=0$ over an arbitrary
field, and strengthened Theorem \ref{sharp3} to show that
finite generation holds if and only if $\rho=0$ \cite{PS}.

{\bf Proof. }
It is classical that, after
blowing up 8 points in $\P^3$ which are a complete intersection
of 3 quadrics, we have an elliptic fibration $f:X\arrow \P^2$;
a reference is Dolgachev-Ortland
\cite[Theorem VI.9]{DO}.
Let us prove that.
We have a net of quadrics
through the given 8 points. So the vector space of sections
of $-\half K_X=2H-E_1-\cdots-E_8$ has dimension 3, and it gives
a rational map $f$ from $X$ to $\P^2$ which is clearly defined
outside $E_1,\ldots,E_8$. Because our three quadrics intersect in 8
distinct points, their intersection must be transverse at each
of the 8 points, and so $f$ is in fact a morphism $X\arrow \P^2$,
with $f^*H\cong -\frac{1}{2}K_X$. It follows that
the canonical bundle
$K_X$ is trivial on all fibers of $f$, and so all smooth
fibers are curves with trivial canonical bundle, that is,
elliptic curves. (Explicitly, all fibers are complete
intersections of two quadrics in $\P^3$.) This construction
also shows that the exceptional divisors $E_1,\ldots,E_8$
are sections of $f$.

By Lemma \ref{minimal},
for any two sections
$S_1$ and $S_2$ of the minimal elliptic fibration $f$ over the general
point of $\P^2$, adding $S_2-S_1$
using the group structure on the general fiber 
is a pseudo-automorphism of $X$ which takes $S_1$ to $S_2$.
Since the exceptional divisor $E_1$ is a section of $f$ and also
a $(-1)$-divisor in Mukai's sense, it follows that
all sections of $f$ over the generic point of $Y$
are $(-1)$-divisors in $X$.

Thus, suppose we can check that the general fiber $E$ 
over $k(\P^2)$ of our elliptic fibration
has infinite Mordell-Weil group. Then $X$ has infinitely many $(-1)$-divisors,
and so the total coordinate ring of $X$ is not finitely generated.
Equivalently, the ring of invariants for the 16-dimensional representation
of $(G_a)^4$ we are considering is not finitely generated.

The Picard group of the general fiber $E$ over the function field $k(\P^2)$
is the quotient of the Picard group of the 3-fold $X$ by
the classes of all irreducible divisors in $X$ which do not map onto $\P^2$.
Since $f^*(H)=-\frac{1}{2}K_X$ and $\Pic(\P^2)=\Z H$, the pullback under
$f:X\arrow \P^2$ of every irreducible divisor in $\P^2$
is a multiple of $-\frac{1}{2}K_X$.
So $\Pic(E)$ is the quotient
of $\Pic(X)\cong \Z^{9}$ by the class of $-\frac{1}{2}K_X$
together with
$r_F$ classes for each irreducible divisor in $\P^2$ whose
inverse image in $X$ has $r_F+1$ irreducible
components, say $\sum_{j=1}^{r_F+1}m_{Fj}D_{Fj}$. Moreover,
the divisors $D_{Fj}$ for all $F$ and $1\leq j\leq r_F$, together
with $-\frac{1}{2}K_X$, are linearly independent in $\Pic(X)_{\Q}$
(even modulo numerical equivalence). This follows from the corresponding
fact about morphisms from a surface to a curve
\cite[Cor.~VIII.4]{Beauville}
by restricting the morphism $X\arrow \P^2$ to the inverse
image of a general
line in $\P^2$. Therefore $\Pic(E)_{\Q}$ has rank
$9-1-\sum r_F$, and the Mordell-Weil group
$\Pic^0(E)_{\Q}$ has rank $7-\sum r_F$.

Let us analyze the subset of $\P^2$ over which $f:X\arrow \P^2$
has reducible fibers. Suppose that some fiber (the intersection
of two quadrics $Q_1$ and $Q_2$ in our net) contains a line $L$.
Let $Q_3$ be a quadric in our net which is not in the pencil
spanned by $Q_1$ and $Q_2$.
Since $Q_1\cap Q_2\cap Q_3=\{p_1,\ldots,p_8\}$ is smooth of dimension 0,
$Q_3$ must intersect
$L$ transversely in 2 points; so $L$ must be the line through
two of the points $p_1,\ldots,p_8$. Since $Q_3$ does not contain $L$,
the pencil spanned by $Q_1$ and $Q_2$ is the unique pencil in our net
whose base locus contains $L$. Thus there are at most $\binom{8}{2}=28$
fibers of $\P^3\dashrightarrow \P^2$ that contain a line.

So, apart from finitely many fibers, every fiber $F$
of $\P^3\dashrightarrow \P^2$ must be a union of irreducible curves
of degree at least 2, with some multiplicities, and with total degree 4.
The only reducible possibility is for $F$ to be a union of two irreducible 
conic curves (curves of degree 2), each with multiplicity 1.
Since every conic in $\P^3$ is contained in a plane, $F$ is contained
(at least as a set) in a reducible quadric $Q_1$, the union of the two planes.
Since the fiber $F$ contains $\{p_1,\ldots,p_8\}$,
$Q_1$ belongs to the net of quadrics through $\{p_1,\ldots,p_8\}$.
It is clear that the intersection of $Q_1$ with every other
quadric in our net is reducible. Equivalently, $Q_1$ corresponds to a line
in $\P^2$ over which every fiber of $X\arrow \P^2$ is reducible.

Thus we find something not at all clear a priori: the elliptic
fibration $X\arrow \P^2$ has reducible fibers over
the union of a finite set
and finitely many lines. 
These lines are in one-to-one correspondence
with the reducible quadrics through $p_1,\ldots,p_8$. Over a general
point of each line, the fiber has exactly two irreducible components.
Thus the number
$\sum r_F$ is equal to the number of reducible quadrics through
$p_1,\ldots,p_8$.

Since a reducible quadric is the union of two planes,
the number $\sum r_F$ is one half the number of coplanar
quadruples among $p_1,\ldots,p_8$. (Recall from the proof
of Lemma \ref{plane} that $p_1,\ldots,p_4$ are coplanar
if and only if $p_5,\ldots,p_8$ are coplanar.) Thus the Mordell-Weil
rank of the elliptic fibration $X\arrow \P^2$ over the generic point
of $\P^2$ is
$\rho=7-\frac{a}{2}$. Therefore $\rho$ is nonnegative.
If $\rho$ is greater than 0, then the Mordell-Weil group is infinite,
$X$ contains infinitely many $(-1)$-divisors, and
the total coordinate
ring of $X$ and the corresponding ring of invariants are infinitely generated.
\qed (Theorem \ref{sharp3})

{\bf Proof of Corollary \ref{four}. }
Let $k$ be any field of characteristic not 2.
By Theorem \ref{mukai}, the ring of invariants in Corollary \ref{four}
is the total
coordinate ring of the blow-up $X$ of $\P^3$ at 8 points, namely
the 8 points given in affine coordinates $w=1$ by $(x,y,z)$ with
$x,y,z\in \{0,1\}$. These 8 points are
the intersection of three (very simple) quadric surfaces,
$x^2=xw$, $y^2=yw$, and $z^2=zw$. To be precise, it might appear
that the above representation of $(G_a)^4$
corresponds to the dual of this arrangement
of 8 points in $\P^3$, but in fact the dual arrangement is projectively
isomorphic to the given one. That is a general property of complete
intersections of 3 quadrics in $\P^3$, by Coble
and Dolgachev-Ortland \cite[Theorem III.3,
Example III.6]{DO}.

Using that the field $k$ has characteristic not 2,
we count that there are exactly 12 coplanar quadruples among these 8
points, corresponding to 6 reducible quadrics in our net:
$x(x-w)$, $y(y-w)$, $z(z-w)$, $(x^2-xw)-(y^2-yw)=(x-y)(x+y-w)$,
$(x-z)(x+z-w)$, and $(y-z)(y+z-w)$. Since $\rho=7-\frac{1}{2}(12)=1$
is greater than 0,
Theorem \ref{sharp3} shows that 
the Mordell-Weil group of the general
fiber $E$ of $X$ over $k(\P^2)$ is infinite, and therefore
the ring of invariants for the given
16-dimensional representation of $(G_a)^4$, over $\F_p$ for $p\geq 3$
or over $\Q$, is not finitely generated.

The Appendix gives a simple example of a 16-dimensional
representation of $(G_a)^4$ over $\F_2$ whose
ring of invariants is not finitely generated. Corollary
\ref{four} is proved. \qed

\section{A generalization of the Kawamata and Morrison conjectures
on Calabi-Yau varieties}

Theorems \ref{sharp2} and \ref{sharp3} relate finite generation
of the total coordinate ring, in some situations,
to finiteness of a certain Mordell-Weil
group. In this section, we make a more general conjecture:
the cone of curves of any variety with semi-ample
anticanonical bundle
should be controlled by a certain group, which may be infinite.
(A line bundle is {\it semi-ample }if some positive multiple
is basepoint-free.)
This would follow from a generalization of the
conjectures of Kawamata and Morrison
on Calabi-Yau fiber spaces to allow klt pairs (Conjecture \ref{domain}).
We prove the conjecture
for smooth projective surfaces with semi-ample anticanonical bundle,
the new case being that of rational elliptic surfaces (Theorem
\ref{rational}).

For a projective morphism $f:X\arrow S$ of normal varieties
with connected fibers, define $N^1(X/S)$ as the real vector
space spanned by Cartier divisors on $X$ modulo numerical equivalence 
on curves on $X$ mapped to a point in $S$.
Define a {\it small $\Q$-factorial modification }(SQM)
of $X$ over $S$ to be a birational map $f:X\dashrightarrow X'$ over $S$,
with $X'$ projective over $S$ and $\Q$-factorial, which is an
isomorphism in codimension one.
A Cartier divisor $D$ on $X$ is called {\it $f$-nef} (resp.\
{\it $f$-movable, $f$-effective}) if $D\cdot C\geq 0$ for every
curve $C$ on $X$ which is mapped to a point in $S$ (resp.,
if $\text{codim}(\text{supp}(\text{coker}(f^*f_*O_X(D)\arrow O_X(D))))\geq 2$,
if $f_*O_X(D)\neq 0$). For an $\R$-divisor $\Delta$ on
a $\Q$-factorial variety $X$,
the pair $(X,\Delta)$ is {\it klt }if, for any resolution
$\pi:\widetilde{X}\arrow X$ with a simple normal crossing
$\R$-divisor $\widetilde{\Delta}$
such that $K_{\widetilde{X}}+\widetilde{\Delta}=\pi^*(K_X+\Delta)$,
the coefficients of $\widetilde{\Delta}$
are less than 1 \cite[Definition 2.34]{KM}.

 The {\it $f$-nef cone }$\A(X/S)$
(resp.\ the {\it closed $f$-movable cone }$\overline{M}(X/S)$)
is the closed convex cone in $N^1(X/S)$ generated by
the numerical classes of $f$-nef divisors (resp.\ $f$-movable divisors).
The {\it $f$-effective cone }$B^e(X/S)$ is the convex cone,
not necessarily closed, generated by $f$-effective Cartier divisors.
We call $A^e(X/S)=\A(X/S)\cap B^e(X/S)$ and
$M^e(X/S)=\overline{M}(X/S)\cap B^e(X/S)$ the {\it $f$-effective
$f$-nef cone }and the {\it $f$-effective $f$-movable cone, }respectively.
Finally, a {\it finite rational polyhedral }cone in $N^1(X/S)$ means the closed
convex cone spanned by a finite set of Cartier divisors on $X$.

\begin{conjecture}
\label{domain}
Let $f:X\arrow S$ be a projective morphism with connected fibers,
$(X,\Delta)$ a $\Q$-factorial klt pair with $\Delta$ effective.
Suppose that $K_X+\Delta$ is numerically trivial over $S$. 
Let $\Aut(X/S)$ and $\PsAut(X/S)$ denote the groups of automorphisms
or pseudo-automorphisms of $X$ over the identity on $S$. Then:

(1) The number of $\Aut(X/S)$-equivalence classes of faces
of the cone $A^e(X/S)$ corresponding to birational contractions
or fiber space structures is finite. Moreover, there exists
a finite rational polyhedral cone $\Pi$ which is a fundamental
domain for the action of $\Aut(X/S)$ on $A^e(X/S)$ in the sense
that

(a) $A^e(X/S)=\cup_{g\in \Aut(X/S)} g_*\Pi$,

(b) $\Int \Pi \cap g_*\Int \Pi=\emptyset$ unless $g_*=1$.

(2) The number of $\PsAut(X/S)$-equivalence classes of chambers
$A^e(X'/S,\alpha)$ on the cone $M^e(X/S)$ corresponding to
marked SQMs $X'\arrow S$ of $X\arrow S$ is finite.
Equivalently, the number of isomorphism classes over $S$
of SQMs of $X$ over $S$ (ignoring the birational
identification with $X$) is finite. Moreover, there exists
a finite rational polyhedral cone $\Pi'$ which is a fundamental
domain for the action of $\PsAut(X/S)$ on $M^e(X/S)$.
\end{conjecture}

Note that Conjecture \ref{domain} would not be true for
Calabi-Yau pairs (pairs $(X,\Delta)$ with $K_X+\Delta\equiv 0$) 
that are log-canonical (or dlt) rather than klt.
Let $X$
be the blow-up of $\P^2$ at 9 very general points. Let $\Delta$
be the proper transform of the unique smooth cubic curve through
the 9 points; then $K_X+\Delta\equiv 0$, and so $(X,\Delta)$ is
a log-canonical Calabi-Yau pair. The surface $X$ contains
infinitely many $(-1)$-curves by Nagata \cite{NagataRat2},
and so the nef cone is not finite polyhedral. But the
automorphism group $\Aut(X)$ is trivial and hence does not
have a finite polyhedral fundamental domain on the nef cone.

The conjecture also fails if we allow the $\R$-divisor $\Delta$ to have
negative coefficients. Let $Y$ be a K3 surface whose cone of curves
is not finite polyhedral, and let $X$ be the blow-up of $Y$ at
a very general point. Then $(X,-E)$ is a klt Calabi-Yau pair,
where $E$ is the exceptional curve. Here the cone of curves of $X$
is not finite polyhedral, but $\Aut(X)$ is trivial.

For $X$ terminal and $\Delta=0$, Conjecture
\ref{domain} is exactly Kawamata's conjecture
on Calabi-Yau fiber spaces,
generalizing Morrison's conjecture on Calabi-Yau
varieties \cite{Kawamata, Morrison1, Morrison2}.
(The group in part (2) can then be described as $\Bir(X/S)$, since
all birational automorphisms of $X$ over $S$ are
pseudo-automorphisms when $K_X$ is numerically trivial over $S$.) Assuming
the base field has characteristic zero, Kawamata's conjecture
is known for $X$ of dimension at most 2, for $X$ of dimension 3
with $S$ of positive dimension, and for a few classes 
of Calabi-Yau 3-folds \cite[Remark 1.13]{Kawamata}.
Some evidence for Conjecture \ref{domain}
in the case $\Delta\neq 0$ is provided by Coble surfaces, smooth projective
rational surfaces such that the linear system $|-K_X|$ is empty
but $|-2K_X|$ is not empty. If there is a smooth divisor
$D$ equivalent to $-2K_X$, then 
$(X,D/2)$ is a klt Calabi-Yau pair.
Then Conjecture \ref{domain} is true for $(X,D/2)$. Indeed, the
double cover of $X$ ramified over $D$ is a K3 surface
\cite[Lemma 6.2]{DZ}, and the conjecture
for $X$ follows from Oguiso-Sakurai \cite[Corollary 1.9]{OS}.

Conjecture \ref{domain} would have strong consequences for varieties
with $-K_X$ semi-ample, as considered in this paper. Let
$X$ be a $\Q$-factorial klt variety with $-K_X$ semi-ample.
For any $m>1$ such that $-mK_X$ is basepoint-free, let $\Delta$
be $1/m$ times a general divisor in the linear system $|-mK_X|$;
then $(X,\Delta)$ is klt \cite[Lemma 5.17]{KM}. Moreover,
$K_X+\Delta\equiv 0$, and so $(X,\Delta)$ is a klt Calabi-Yau pair.
Thus Conjecture \ref{domain} implies that conclusions (1) and (2) hold
whenever $-K_X$ is semi-ample. If $-K_X$ is ample, then the
conjecture is known in characteristic zero:
the nef cone is finite rational polyhedral by the cone theorem
\cite{KM}, and the movable cone is finite rational polyhedral,
partitioned into the nef cones of the finitely many SQMs
of $X$,
by Birkar-Cascini-Hacon-M\textsuperscript{c}Kernan \cite{BCHM}.

We can now prove Conjecture \ref{domain} when $X$ is a rational elliptic
surface; this completes the proof for all smooth projective
surfaces with $-K_X$ semi-ample. For rational elliptic surfaces
with no multiple fibers and Mordell-Weil rank 8 (the maximum possible),
Theorem \ref{rational}
was already known, by Grassi and Morrison \cite[Theorem 2.3]{GM}.

\begin{theorem}
\label{rational}
Let $X$ be a rational elliptic (or quasi-elliptic)
surface over an algebraically
closed field. That is, $X$ is a smooth
projective surface with $-K_X$ semi-ample such that
$T:=\Proj R(X,-K_X)$ has dimension 1. Then
the number of $\Aut(X/T)$-equivalence classes of faces
of the nef effective cone $A^e(X)$ corresponding to birational contractions
or fiber space structures is finite. (In this case the nef effective
cone is closed, thus equal to the nef cone.)
Moreover, there exists
a finite rational polyhedral cone $\Pi$ which is a fundamental
domain for the action of $\Aut(X/T)$ on $A^e(X)$.

Also, the following are equivalent:

(1) The total coordinate ring of $X$ is finitely generated.

(2) The nef cone of $X$ is finite rational polyhedral.

(3) The Mordell-Weil group, $\Pic^0$ of the generic fiber of $X$
over $T$, is finite. (The Mordell-Weil group is automatically
a finite-index subgroup of $\Aut(X/T)$.)
\end{theorem}

Note that we only have to consider statement (1) in Conjecture \ref{domain};
statement (2) in the conjecture is vacuous for surfaces, because
every movable divisor on a surface is nef. (Or, related to that:
minimal models of surfaces are unique.)

{\bf Proof. }The basic point is to show that the group
$\Aut(X/T)$ has only finitely many orbits on the set
of $(-1)$-curves in $X$.

The generic fiber of $X\arrow T$ is a curve
$X_{\eta}$ of genus 1 over the function field of $T$. The class
in $\Pic (X)$ of a general fiber of $X\arrow T$ is $-mK_X$
for some positive integer $m$. Here $-K_X$ is effective by Riemann-Roch,
and so (if $m>1$) there is a multiple fiber of $X\arrow T$, a curve
in the class of $-K_X$. (Here $m=1$ if and only if $X\arrow T$
has a section \cite[Chapter 5, \S 6]{CD}.)
The irreducible components
of the reducible fibers of $X\arrow T$, if any, are exactly the
$(-2)$-curves on $X$ (smooth curves $C$ of genus 0
with $K_X\cdot C=0$). The Picard group $\Pic (X_{\eta})$ is the quotient
of $\Pic (X)$ by $-K_X$ together with
all the $(-2)$-curves. The degree of a line bundle on $X$
on a general fiber of $X\arrow T$ is given by the intersection
number with $-mK_X$, and so the Mordell-Weil group $G:=\Pic^0(X_{\eta})$
is the subquotient of $\Pic (X)$ given by
$$G=(K_X)^{\perp}/(-K_X, (-2)\text{-curves}).$$

An element $x$ of the group $G$
acts by a translation on the curve $X_{\eta}$ of genus 1
and hence by an automorphism $\varphi_x$
on $X$, by Lemma \ref{minimal}.
This gives an action of $G$ on $\Pic (X)$.
We know how translation by an element $x$ of $\Pic^0(X_{\eta})$
acts on $\Pic (X_{\eta})$: by $\varphi_x(y)=
y+\deg(y)x$. Since $-mK_X\in \Pic (X)$ is the class of a general
fiber of $X\arrow T$, this means that $\varphi_x$ acts on $\Pic(X)$
by $\varphi_x(y)=y-(mK_X\cdot y) x \pmod{-K_X,
(-2)\text{-curves}}$. Using that the action of $G$ preserves
the intersection product, we compute the action of certain
elements of $G$ on $\Pic(X)$ by:
$$\varphi_x(y)=y -(mK_X\cdot y )x + \big[ y\cdot
(-x+(1/2)(x\cdot x) mK_X)\big] (-mK_X)$$
for all $x\in K_X^{\perp}$ with $x\cdot C_i=0$ for all
$(-2)$-curves $C_i$, and all $y\in \Pic(X)$.

Now let $E_1$ and $E_2$ be any two $(-1)$-curves on $X$
such that $E_1\cdot C_i=E_2\cdot C_i$ for all $(-2)$-curves
$C_i$ and such that $E_1\equiv E_2\pmod{m\Pic(X)}$.
Let $x=(E_2-E_1)/m\in \Pic(X)$. Then $x$ is in $K_X^{\perp}$,
we have
$x\cdot C_i=0$ for all $(-2)$-curves $C_i$, and
\begin{align*}
\varphi_x(E_1) &= E_1-(mK_X\cdot E_1) x
+ \big[ E_1\cdot (-mx+(1/2)(mx\cdot mx) K_X)\big] (-K_X)\\
 &= E_2.
\end{align*}

Every $(-1)$-curve $E$ has $E\cdot(-mK_X)=m$, which says that
$E$ is a multisection of degree $m$ of $X\arrow T$.
Therefore $0\leq E\cdot C_i\leq m$ for each $(-2)$-curve $C_i$.
So the $(-1)$-curves $E$ are divided into finitely many
classes according to the intersection numbers of $E$
with all $(-2)$-curves
and the class of $E$ in $\Pic(X)/m$. By the previous paragraph,
the $(-1)$-curves on $X$ fall into finitely many
orbits under the action of $G$. A fortiori, $\Aut(X/T)$ has finitely
many orbits on the set of $(-1)$-curves.

We now describe all the extremal rays of the cone of curves $\NE(X)$.
By the cone theorem \cite[2.1.1, 2.3.3]{KK},
all the $K_X$-negative extremal rays of $\NE(X)$ are spanned by $(-1)$-curves.
Since $-K_X$ is nef, it remains to describe $\NE(X)\cap K_X^{\perp}$.
It is a general fact, for any smooth projective surface $X$,
that any extremal ray $\R^{>0}\cdot x$ of $\NE(X)$ with
$x^2<0$ is spanned by a curve $C$ with $C^2<0$, and that such
an extremal ray is isolated. Both statements follow from
the fact that any two distinct curves $C$ and $D$ on $X$ with
$C^2<0$ and $D^2<0$ have $C\cdot D\geq 0$ and hence are
``far'' from each other. 

For the rational elliptic
surface $X$, the Hodge index theorem gives that
the intersection pairing on $K_X^{\perp}$ is negative
semidefinite, with $x^2=0$ only on the line spanned by $-K_X$.
Therefore, using the fact in the previous paragraph, all extremal rays
of $\NE(X)\cap K_X^{\perp}$ are spanned by either $-K_X$
or a $(-2)$-curve. There are only finitely many $(-2)$-curves
(the irreducible components of reducible fibers of $X\arrow T$),
and so the cone $\NE(X)\cap K_X^{\perp}$ is finite rational
polyhedral. We conclude that every extremal ray
of $\NE(X)$ is spanned by either
a $(-1)$-curve or one of 
the finitely many $(-2)$-curves (or $-K_X$,
if $X$ contains no $(-2)$-curve).

Since $\Aut(X/T)$ has only finitely many orbits on $(-1)$-curves,
it follows that $\Aut(X/T)$ has only finitely many orbits
on the extremal rays of $\NE(X)$.
Moreover, if $X$ contains infinitely many $(-1)$-curves,
then the only possible limit ray of $(-1)$-rays is
$\R^{>0}(-K_X)$. Indeed, every $(-1)$-curve $E$ has $-K_X\cdot E=1$ and
$E^2=-1$,
while there are only finitely many $(-1)$-curves $E$ with
any given degree $H\cdot E$ (a natural number), where $H$
is a fixed ample divisor on $X$. Any limit ray
$\R^{>0}x$ of $(-1)$-rays has $0<H\cdot x<\infty$, and so it must
have $(-K_X)\cdot x=0$ and $x^2=0$. Since the intersection form
is negative semidefinite on $(K_X)^{\perp}$ with kernel
spanned by $-K_X$, it follows that $x$ is a multiple of $-K_X$,
as claimed.

We can deduce that the nef cone $\A(X)$ is finite rational
polyhedral near any point $x$ in $\A(X)$ not in the ray
$\R^{\geq 0}(-K_X)$. First, such a point $x$ has $x^2\geq 0$
and also $(-K_X)\cdot x\geq 0$, since $x$ and $-K_X$ are nef. If
$(-K_X)\cdot x$ were equal to zero, these properties would imply
that $x$ is a multiple of $-K_X$; hence we must have $(-K_X)\cdot x>0$.
As a result, there is a neighborhood $V$ of $x$ and a
neighborhood $W$ of $-K_X$ such that $V\cdot W>0$. Since
the only possible limit ray of $(-1)$-rays is $\R^{>0}(-K_X)$,
almost all (all but finitely many) $(-1)$-curves are in the cone $\R^{>0} W$.
So almost all $(-1)$-curves have positive intersection
with the neighborhood $V$ of $x$. Since almost all
extremal rays of $\NE(X)$ are spanned by $(-1)$-curves,
we conclude that the nef cone $\A(X)$ is finite rational polyhedral
near $x$, as claimed.

In particular, for each $(-1)$-curve $E$, the face $\A(X)\cap E^{\perp}$
of the nef cone is finite rational polyhedral, since it does not
contain $-K_X$. So the cone $\Pi_E$ spanned by $-K_X$
and $\A(X)\cap E^{\perp}$ is finite rational polyhedral.

Let $x$ be any nef $\R$-divisor on $X$.
Let $c$ be the maximum real number such that $y:=x+cK_X$
is nef. Then $x$ and $y$ have the same degree on all $(-2)$-curves,
and so there must be some $(-1)$-curve $E$ with $y\in E^{\perp}$.
Therefore $x$ is in the cone $\Pi_E$. That is, the nef cone
$\A(X)$ is the union of the finite rational polyhedral cones
$\Pi_E$. Moreover, for two distinct $(-1)$-curves $E$ and $F$
on $X$, the intersection $\Int \Pi_E\cap \Int \Pi_F$ is empty,
since $E$ and $F$ are linearly independent in $N^1(X)$. Since
$\Aut(X/T)$ has only finitely many orbits on the set of
cones $\Pi_E$,
there is a finite rational polyhedral cone $\Pi$ which is
a fundamental domain for the action of $\Aut(X/T)$ on $A^e(X)$.
(We can construct such a domain using the Dirichlet construction
for the hyperbolic metric on $\{x\in N^1(X): x^2>0, H\cdot x>0\}/
\R^{>0}$, restricted to the nef cone. That is, given a point
$x$ in the interior of the nef cone, we define a fundamental domain
as the set of points of the
nef cone whose distance to $x$, in the nef cone modulo scalars,
is at most their
distance to any
other point in the $\Aut(X/T)$-orbit of $x$.)

Any rational point $x$ in the nef cone $\A(X)$ is effective. This
is clear from the Riemann-Roch theorem $\chi(X,nL)=\chi(X,O)
+((nL)^2-K_X\cdot nL)/2$, since for a nef divisor $x$ on our
surface $X$, either $x^2>0$, $x^2=0$ and $(-K_X)\cdot x>0$, or
$x$ is a multiple of $-K_X$. Since
the cone $\Pi_E$ is finite rational polyhedral for
each $(-1)$-curve $E$, 
$\Pi_E$ is contained in the nef effective cone $A^e(X)$.
Since the whole nef cone $\A(X)$ is the union of the cones
$\Pi_E$, the nef effective cone $A^e(X)$ is equal to its closure,
the nef cone $\A(X)$. This proves another statement
of Theorem \ref{rational}.

Next, we have to show that the number of $\Aut(X/T)$-orbits
of faces of the nef effective cone $A^e(X)$
corresponding to birational contractions or fiber
space structures is finite. Dually, it suffices to show
that there are only finitely many $\Aut(X/T)$-orbits of
faces in the cone of curves $\NE(X)$. Since there are only
finitely many orbits of $(-1)$-curves, and almost all
extremal rays of $\NE(X)$ are spanned by $(-1)$-curves, it
suffices to show that each $(-1)$-curve $E$ lies on only finitely
many faces of $\NE(X)$. This follows from the dual statement
that the cone $\A(X)\cap E^{\perp}$ is finite rational polyhedral,
which we have proved.

Finally, let us prove the equivalence of statements (1), (2) and (3)
in Theorem \ref{rational}.
That (1) implies (2) is easy; it is also part of Hu and Keel's
characterization of varieties with finitely generated total coordinate ring
\cite[Proposition 2.9]{HK}. Also, (2) easily implies (3),
as follows. Suppose that the Mordell-Weil group $\Pic^0(X_{\eta})$
is infinite. We know that $X$ contains at least one $(-1)$-curve
since $X$ is a rational elliptic surface. The Mordell-Weil
group $\Pic^0(X_{\eta})$ acts on $X$, and only a finite subgroup
of it can map any given $(-1)$-curve into itself, since a $(-1)$-curve
is a multisection of $X\arrow T$ of degree a positive integer $m$.
Thus, if $\Pic^0(X_{\eta})$ is infinite, then $X$ has infinitely
many $(-1)$-curves, and so $\NE(X)$ is not finite rational polyhedral.
So the dual cone $\A(X)$ is not finite rational polyhedral.
That is, (2) implies (3).

Conversely, suppose (3), that $\Pic^0(X_{\eta})$ is finite. This
is a finite-index subgroup of $\Aut(X/T)$, and so $\Aut(X/T)$
is finite. By the earlier parts of this theorem, it follows
that the nef effective cone $A^e(X)$ is finite rational
polyhedral. Statement (2) is proved. It remains to prove (1),
that the total coordinate ring of $X$ is finitely generated.
By Hu-Keel's theorem \cite[Proposition 2.9]{HK}, it suffices
to show that every codimension-1 face of $\NE(X)$ can
be contracted. This follows from the cone theorem (including
the contraction theorem), by the same argument as in the proof
of Theorem \ref{sharp2}. \qed

\section{Appendix: Infinitely generated rings of invariants over
$\F_2$ or $\F_3$}

In this section we complete the proof of Corollaries \ref{six}
and \ref{four} by giving
simple examples of 18-dimensional
representations of $(G_a)^6$ over $\F_2$ or $\F_3$,
and 16-dimensional representations of $(G_a)^4$ over $\F_2$,
whose rings of invariants are not finitely generated.

Take the 9 points in $\P^2(\F_4)$ given in affine coordinates
by $(x,y)$ where $x$ and $y$ run through the nonzero elements of a field
$\F_4$ of order 4.
These points are the intersection of the two cubics $x^3=z^3$
and $y^3=z^3$, each consisting of three lines through a point in $\P^2$.
There are 9 collinear triples among these points, and
3 partitions of them into disjoint collinear triples. (These partitions
correspond to three cubics in the pencil which are unions of three lines,
$x^3+z^3$, $y^3+z^3$, and $(x^3+z^3)+(y^3+z^3)=x^3+y^3$.) Therefore
$\rho=8-9+3=2$ is greater than 0. By Theorem \ref{sharp2},
the blow-up $X$ of $\P^2$ at these 9 points has infinitely many
$(-1)$-curves. (More precisely, this count shows that
the Mordell-Weil rank is 2.)
Thus we get an 18-dimensional representation
of $(G_a)^6$ over $\F_4$ such that the ring of invariants is not
finitely generated.

Since this pencil of cubics is defined over $\F_2$, even though the
9 individual points are not, this pencil gives an 18-dimensional
representation of $(G_a)^6$ over $\F_2$ such that the ring of invariants
is not finitely generated, as explained in section \ref{allfields}.
We will write out this representation explicitly, although the
reader might be satisfied to know that it can be done.
Let $\zeta$
be a primitive cube root of unity in $\F_4$.
Our representation
of $(G_a)^6$ on $A^{18}:=A^9_1\oplus A^9_2$
over $\F_4$ can be defined by 6 commuting
linear maps $1+A_i$, $1\leq i\leq 6$, where $A_i\in \Hom(A^9_2,A^9_1)$
are given by the following diagonal matrices:
\begin{align*}
\diag(1,\zeta,\zeta^2,0,0,0,0,0,0),& \quad\diag(1,0,0,\zeta,\zeta^2,0,0,0,0),\\
\diag(1,0,0,0,0,\zeta,\zeta^2,0,0),& \quad 
\diag(0,1,1,1,1,1,1,0,0),\\
\diag(0,0,0,0,0,1,1,1,1),& \quad\diag(1,1,1,0,0,0,0,
\zeta,\zeta^2).
\end{align*}
These 6 vectors are chosen as a basis for the kernel of the
linear map $A^9\arrow
A^3$ corresponding to our 9 points in $\P^2$ over $\F_4$,
in the order $[1,1,1]$, $[1,\zeta,1]$, $[1,\zeta^2,1]$,
$[\zeta,1,1]$, $[\zeta^2,1,1]$, $[\zeta,\zeta,1]$,
$[\zeta^2,\zeta^2,1]$, $[\zeta,\zeta^2,1]$,
$[\zeta^2,\zeta,1]$.
We can define a form of this representation over $\F_2$
by 6 commuting linear maps $1+A_i$, $1\leq i\leq 6$, where
$A_i\in \Hom(A^9_2,A^9_1)$ are given by the following matrices:
\begin{align*}
\diag(1,\begin{pmatrix}0 & 1\\ 1 & 1 \end{pmatrix},0,0,0,0,0,0), &\quad
\diag(1,0,0,\begin{pmatrix}0 & 1\\ 1 & 1 \end{pmatrix},0,0,0,0),\\
\diag(1,0,0,0,0,\begin{pmatrix}0 & 1\\ 1 & 1 \end{pmatrix},0,0), &\quad
\diag(0,1,1,1,1,1,1,0,0),\\
\diag(0,0,0,0,0,1,1,1,1), &\quad
\diag(1,1,1,0,0,0,0,\begin{pmatrix}0 & 1\\ 1 & 1 \end{pmatrix}).
\end{align*}
This is an 18-dimensional representation of $(G_a)^6$ over $\F_2$
with infinitely generated ring of invariants.

A similar example works over $\F_3$. Take the
9 points of $\P^2(\F_9)$ given in affine coordinates by $(x,y)$
where $x$ and $y$ run through the set of fourth roots of unity in $\F_9$
other than 1. These points are the intersection of the two cubics
$x^3+x^2z+xz^2+z^3$ and $y^3+y^2z+yz^2+z^3$. There are 7
lines through these points, and 2 partitions of them into
disjoint collinear triples. (The two partitions correspond
to two cubics in the pencil which are unions of three lines,
$x^3+x^2z+xz^2+z^3$ and $y^3+y^2z+yz^2+z^3$, and the other
collinear triple corresponds to another reducible cubic
in the pencil,
$(x^3+x^2z+xz^2+z^3)-(y^3+y^2z+yz^2+z^3)=(x-y)(x^2+xy+y^2+xz+yz+z^2)$.)
So $\rho=8-7+2=3$. By Theorem \ref{sharp2},
the blow-up $X$ of $\P^2$ at these 9 points has infinitely
many $(-1)$-curves. This gives an 18-dimensional representation
of $(G_a)^6$ over $\F_9$ whose ring of invariants is not
finitely generated. 

Since this pencil of cubics is defined over $\F_3$, we can
define a form of the above representation over $\F_3$.
Explicitly, let $i$ be a square root of $-1$ in $\F_9$. Then our
representation
of $(G_a)^6$ on $A^{18}:=A^9_1\oplus A^9_2$
over $\F_9$ can be defined by 6 commuting
linear maps $1+A_j$, $1\leq j\leq 6$, where $A_j\in \Hom(A^9_2,A^9_1)$
are given by the following diagonal matrices:
\begin{align*}
\diag(1,1+i,1-i,0,0,0,0,0,0),& \quad \diag(1,0,0,1+i,1-i,0,0,0,0),\\
\diag(1,0,0,0,0,1+i,1-i,0,0),& \quad
\diag(1,1,1,1,1,-1,-1,0,0),\\
\diag(0,0,0,0,0,1,1,-1,-1),& \quad \diag(0,-1,-1,1,1,-1,-1,1+i,1-i).
\end{align*}
These 6 vectors are chosen as a basis for the kernel of the
linear map $A^9\arrow
A^3$ corresponding to our 9 points in $\P^2$ over $\F_9$,
in the order $[-1,-1,1]$, $[-1,i,1]$, $[-1,-i,1]$, $[i,-1,1]$,
$[-i,-1,1]$,
$[i,i,1]$, $[-i,-i,1]$, $[i,-i,1]$, $[-i,i,1]$.
We can define a form of this representation over $\F_3$
by 6 commuting linear maps $1+A_j$, $1\leq j\leq 6$, where
$A_j\in \Hom(A^9_2,A^9_1)$ are given by the following matrices:
\begin{align*}
\diag(1,\begin{pmatrix}0 & 1\\ 1 & -1 \end{pmatrix},0,0,0,0,0,0), &\quad
\diag(1,0,0,\begin{pmatrix}0 & 1\\ 1 & -1 \end{pmatrix},0,0,0,0),\\
\diag(1,0,0,0,0,\begin{pmatrix}0 & 1\\ 1 & -1 \end{pmatrix},0,0), &\quad
\diag(1,1,1,1,1,-1,-1,0,0),\\
\diag(0,0,0,0,0,1,1,-1,-1), &\quad
\diag(0,-1,-1,1,1,-1,-1,\begin{pmatrix}0 & 1\\ 1 & -1 \end{pmatrix}).
\end{align*}
This is an 18-dimensional representation of $(G_a)^6$ over $\F_3$
with infinitely generated ring of invariants. Corollary
\ref{six} is proved.

We now finish the proof of Corollary \ref{four} by exhibiting
a simple example of a 16-dimensional
representation of $(G_a)^4$ over $\F_2$ whose
ring of invariants is not finitely generated.

Let $\zeta$
be a primitive cube root of unity in $\F_4$, and consider the 8 points
in $\P^3(\F_4)$ defined in affine coordinates $w=1$ by $(x,y,z)$ where
$z$ is $\zeta$ or $\zeta^2$, $y$ is 0 or $z$, and $x$ is 0 or 1. 
These points are the complete intersection of the three quadrics
$x^2+xw=0$, $y^2+yz=0$, and $z^2+zw+w^2=0$. Let $X$ be the blow-up
of $\P^3$ at these 8 points. We have an elliptic fibration
$X\arrow \P^2$. We count that there are exactly 10
coplanar quadruples among the 8 points,
corresponding to 5 reducible quadrics in our net:
$x(x+w)$, $y(y+z)$, $z^2+zw+w^2=(z+\zeta w)(z+\zeta^2 w)$,
$y^2+yz+z^2+zw+w^2=(y+w+\zeta z)(y+w+\zeta^2 z)$,
and $x^2+xw+z^2+zw+w^2=(x+z+\zeta w)(x+z+\zeta^2 w)$. Since
$\rho=7-\frac{1}{2}(10)=2$ is greater than 0,
Theorem \ref{sharp3} shows that 
the Mordell-Weil group of the general
fiber $E$ of $X$ over $k(\P^2)$ is infinite (in fact of rank 2),
and therefore
the ring of invariants for the given
16-dimensional representation of $(G_a)^4$ over $\F_4$
is not finitely generated.

Since this net of quadrics is defined over $\F_2$, even though the
8 individual points are not, we can define a form of the above
representation over $\F_2$. Our representation
of $(G_a)^4$ on $A^{16}:=A^8_1\oplus A^8_2$
over $\F_4$ can be defined by 4 commuting
linear maps $1+A_i$, $1\leq i\leq 4$, where $A_i\in \Hom(A^8_2,A^8_1)$
are given by the following diagonal matrices:
\begin{align*}
\diag(0,0,1,1,0,0,1,1),& \quad\diag(0,0,0,0,\zeta,\zeta^2,\zeta,\zeta^2),\\
\diag(\zeta,\zeta^2,\zeta,\zeta^2,\zeta,\zeta^2,\zeta,\zeta^2),& \quad 
\diag(1,1,1,1,1,1,1,1).\\
\end{align*}
These 4 vectors are chosen as a basis for the kernel of the
linear map $A^8\arrow
A^4$ corresponding to our 8 points in $\P^2$ over $\F_4$,
in the order $[0,0,\zeta,1]$, $[0,0,\zeta^2,1]$, $[1,0,\zeta,1]$,
$[1,0,\zeta^2,1]$, $[0,\zeta,\zeta,1]$, $[0,\zeta^2,\zeta^2,1]$,
$[1,\zeta,\zeta,1]$, $[1,\zeta^2,\zeta^2,1]$.
We can define a form of this representation over $\F_2$
by 4 commuting linear maps $1+A_i$, $1\leq i\leq 4$, where
$A_i\in \Hom(A^8_2,A^8_1)$ are given by the following matrices:
\begin{align*}
\diag(0,0,1,1,0,0,1,1), &\quad 
\diag(0,0,0,0,\begin{pmatrix}0 & 1\\ 1 & 1 \end{pmatrix},\begin{pmatrix}0 & 1\\ 1 & 1 \end{pmatrix}), \\
\diag(\begin{pmatrix}0 & 1\\ 1 & 1 \end{pmatrix},\begin{pmatrix}0 & 1\\ 1 & 1 \end{pmatrix},\begin{pmatrix}0 & 1\\ 1 & 1 \end{pmatrix},\begin{pmatrix}0 & 1\\ 1 & 1 \end{pmatrix}) , &\quad \diag(1,1,1,1,1,1,1,1).
\end{align*}
This is a 16-dimensional representation of $(G_a)^4$ over $\F_2$
with infinitely generated ring of invariants. That completes
the proof of Corollary \ref{four}.

% Omit these bibliography lines if there's no bibliography.

\small \sc DPMMS, Wilberforce Road,
Cambridge CB3 0WB, England

b.totaro@dpmms.cam.ac.uk
\end{document}